\documentclass[11pt]{book} 
\usepackage{amsmath,amssymb}
\def\bbr{{\mathbb R}}
\newtheorem{lemma}{Lemma}[section]
\newtheorem{prop}{Proposition}[section]
\newtheorem{coro}[prop]{Corollary}
\newtheorem{thm}{Theorem}[section]
\newtheorem{cor}[thm]{Corollary}
\newtheorem{remn}{Remark}

\newenvironment{pf}{\noindent{\it Proof}\rm.}{\hfill q.e.d.}
\numberwithin{equation}{section}
\begin{document}
\chapter*{\Large\rm
COMMUTATORS, EIGENVALUE GAPS, AND MEAN CURVATURE 
IN THE THEORY OF SCHR\"ODINGER OPERATORS\\[1em]
\centerline{\large EVANS M.  HARRELL II\footnotemark}}
\footnotetext{harrell@math.gatech.edu, 
School of Mathematics,
Georgia Institute of Technology, Atlanta, GA 30332-0160, USA\\
\noindent\copyright 2002, 2003 by the author.
Reproduction of this article by any means, in its entirety including 
this notice, is permitted for non-commercial purposes.  
}
\centerline{\bf Abstract}
{\narrower\small
Commutator relations are used to investigate the spectra 
of Schr\"odinger Hamiltonians, $H = -\Delta + V\left({x}\right),$
acting on functions of a smooth, compact $d$-dimensional manifold $M$ 
immersed in $\bbr^{\nu}, \nu \geq d+1$.   
Here $\Delta$ denotes the Laplace-Beltrami
operator, and the real-valued potential--energy function $V(x)$
acts by multiplication.   The manifold $M$ may be complete or it may 
have a boundary, in which case Dirichlet boundary
conditions are imposed.

It is found that the mean curvature of a manifold poses tight constraints 
on the spectrum of $H$.
Further, a special algebraic r\^ole is found to be
played by a Schr\"odinger operator with potential proportional to the
square of the mean curvature: 
$$H_{g} := -\Delta +  g h^2,$$
where $\nu = d+1$, $g$ is a real parameter, and 
$$h := \sum\limits_{j = 1}^{d} {\kappa_j},$$ 
with 
$\{\kappa_j\}$, $j = 1, \dots, d$ denoting the principal curvatures of $M$.
For instance, by Theorem~\ref{thm3.1} and Corollary~\ref{cor4.5}, 
each eigenvalue gap of an arbitrary Schr\"odinger
operator is bounded above by an expression using 
$H_{1/4}$.  The ``isoperimetric" parts of these theorems state that these bounds 
are sharp
for the fundamental eigenvalue gap and for infinitely many other eigenvalue gaps.
\par}

\section{Introduction}

Eigenvalue gaps for Laplacians and Schr\"odinger 
Hamiltonians, 
$$H = -\Delta + V\left({x}\right),$$
especially the fundamental gap, 
	$$\Gamma := \lambda_2 - \lambda_1,$$
are connected to the analytic and geometric
details of $H$ in many familiar ways, yet
the most elementary connections between eigenvalue gaps and
operators are algebraic, through calculations of commutators.
Commutators and their relation to eigenvalue gaps 
and inverse spectral problems
have been studied in \cite{Ash02,AsHe03,Har88,Har93,HaSt97,Her99,LePa02},
which partly inspired this work.
Here $H$ will be defined on a hypersurface, and  commutators 
will be used to connect  eigenvalue
gaps and certain other functions defined on the spectrum 
$\sigma(H)$ 
to the mean curvature of the manifold $M$, which is found
to pose tight constraints 
on the spectrum.

Furthermore, Schr\"odinger operators 
of specific kinds depending explicitly on curvature,
\begin{equation}\label{eq1.1}
-{\nabla }^{2} + q\left({\kappa }\right),\end{equation} 
defined on a $d$-dimensional submanifold
$\bbr^{\nu}$ having principal curvatures
$\kappa := \{\kappa_j\}$, $j = 1, \dots, d$,
arise in the mathematics of nanophysics and thin structures.
Coincidentally, simplifications are achieved when the potential 
energy is of this form.  
Low-lying 
eigenvalues, $\lambda_1 < \lambda_2 \leq \cdots$, 
of (1.1) have been studied, and
sharp estimates have been discovered relating
to the geometry of $M$
in  \cite{Exn90,AsEx90,DuEx95,Har96,HaLo98,ExHaLo99}.
For certain operators of the type
\eqref{eq1.1} with quadratic functions $q$,
the estimates of the gap
$\Gamma$ in this article will likewise be sharp.

Section 2 will be devoted to the fundamental eigenvalue gap in
the one--dimensional situation, where $M$ is a space curve.
In Section 3, a similar analysis will be carried out when 
$d > 1$, while in Section 4 an alternative algebraic
approach is used to extend these results to all eigenvalue gaps.
The results of each section are both generalized and 
largely subsumed
by the later sections, but as alternative,
simpler methods are available for the simpler situations, 
it is felt useful to present multiple points of view.

%
%
\section*{Notation}
\begin{itemize}
\item[]
$\langle{f,g}\rangle$ denotes the standard inner product on $L^2(M)$.
\item[]
$[A,B] := AB - BA$  denotes the commutator of operators $A$ and $B$
\item[]
${\bf t},{\bf n},$ and ${\bf b}$ are the moving trihedron of unit vectors for 
a space curve in $\bbr^3$.   
The corresponding curvature and torsion are $\kappa$ and $\tau$.  
\item[]
When the dimension $d \geq 2$, 
$M$ will denote a compact, smooth, connected manifold of dimension $d$ immersed in
$\bbr^{d+1}$, with principal curvatures $\{\kappa_j\}$, $j = 1, \dots, d$.  The corresponding
normalized eigenvectors of the shape operator will be denoted ${\bf t}_j$, and
$h := \sum\limits_{j = 1}^{d} {\kappa_j}$.
\item[]
The notation $u_k, \lambda_k$, $k = 1, \dots$,  
will be used for the normalized eigenfunctions 
and associated eigenvalues, in increasing order, of $H$,
an operator on a Hilbert space
${\cal H}$.  When $H$ is a Schr\"odinger Hamiltonian, it will be
assumed without loss of generality that $\{u_k\}$ are real-valued, and 
that $u_1$ is nonnegative.  
\end{itemize}

Eigenvalue gaps are connected with commutators at a rather 
general level.  Assume
for now only that $H$ and $G$ 
are symmetric operators on a Hilbert space, 
that $\sigma(H)$ contains discrete eigenvalues, and that certain 
products of operators are well-defined.  (Assumptions are 
spelled out below in Lemma~\ref{lem1.1}.)  
The most basic 
relationship of this kind is the elementary commutator gap formula,
\begin{equation}\label{eq1.2}
\left\langle{u_j, [H,G] u_k}\right\rangle = 
(\lambda_j -\lambda_k) \left\langle{u_j, G u_k}\right\rangle.
\end{equation} 

It turns out that not only the $\left[H,G\right]$ but also 
higher-order commutators such as  $\left[{G,\left[{H,G}\right]}\right]$
are related to  eigenvalue gaps for second-order differential operators.
The basis for that analysis can be recapitulated as follows.

For the first commutator of $H$ with $G,$
since $[H,G] u_k = (H - \lambda_k)G u_k$, 
it follows that
\begin{equation}\label{eq1.3}
\left\|[H,G]{u}_{k}\right\|^{2}=\left\langle G{u}_{k},\left(H-\lambda 
_{k}\right)^{2}G{u}_{k}\right\rangle ,
\end{equation} 
and more generally
\begin{equation}\label{eq1.4}
\left\langle{[H,G]{u}_{j},[H,G]{u}_{k}}\right\rangle= 
\left\langle{G{u}_{j},\left({H-{\lambda }_{j}}\right)\left({H-{\lambda 
}_{k}}\right)G{u}_{k}}\right\rangle.\end{equation} 
For the second commutator, a similar short calculation shows that, formally,
\begin{equation}
\left\langle{{u}_{j}\mid \left[{G,\left[{H,G}\right]}\right]{u}_{k}}
\right\rangle 
= \left\langle{{Gu}_{j}\mid \left({2 H-{\lambda }_{j}-{\lambda 
}_{k}}\right){Gu}_{k}}\right\rangle.\end{equation} 
In particular,
\begin{equation}\label{eq1.5}
\left\langle{{u}_{j}\mid \left[{G,\left[{H,G}\right]}\right]{u}_{j}}
\right\rangle 
= 2 \left\langle{{Gu}_{j}\mid \left({H-{\lambda }_{j}}\right){Gu}_{j}}\right\rangle.
\end{equation}

\begin{lemma}\label{lem1.1} 
Let $H$ be a positive self-adjoint operator with discrete 
eigenvalues $\lambda_1$ and $\lambda_2$.  
Let $P$ denote the orthogonal projection 
onto $u_1$, and suppose $G$ is a self-adjoint
operator such that the products $GP$, $G^2P$, $HG^2P$, $H^2GP$,
 and $GHGP$ are 
defined.  Then the fundamental gap $\Gamma := \Gamma(H)$ satisfies
\begin{equation}\label{eq1.6}
\Gamma \left\langle{{u}_{1}, \left[{G,[H,G]}\right]{u}_{1}}\right\rangle 
\le  2 {\left\|[H,G]{u}_{1}\right\| }^{2} .\end{equation} 
\end{lemma}

\begin{remn} 
\rm This is actually a special case of 
Theorem 2.1 of \cite{Har93}.  A simplified proof
is presented here for convenience.
\end{remn}

\begin{pf} 
The assumptions on the products justify formulae 
\eqref{eq1.3}--\eqref{eq1.5}.  With $k=j=1,$
the inequality is thus equivalent to
$$\left\langle{G{u}_{1}\mid \left({H-{\lambda 
}_{1}}\right)\left({{\lambda }_{2}-{\lambda 
}_{1}}\right)G{u}_{1}}\right\rangle \le
\left\langle{G{u}_{1}\mid {\left({H-{\lambda 
}_{1}}\right)}^{2}G{u}_{1}}\right\rangle, $$
which follows from the spectral functional calculus 
(e.g., see \cite{Dav95, EdEv87, ReSi80}),
since for $\mu \in \sigma(H)$, 
$\left(\mu-\lambda_1)\right)\left(\lambda_2-\lambda_1)\right) 
\leq \left(\mu-\lambda_1)\right)^2$.
\end{pf}

\medskip
Whereas identities \eqref{eq1.3}--\eqref{eq1.5} and 
Lemma~\ref{lem1.1} are general facts about operator algebra,
if $H$ is a differential operator, then its
commutators satisfy further algebraic relations.
For instance, if  $H = -\Delta$ is a Laplace operator on 
a Euclidean set $\Omega$, then with
the choice $G = x_k$:
\begin{equation}\label{eq1.7}
\sum\limits_{k=1}^{\nu } {\left[{H, {x}_{k}}\right]}^*
\left[{H, {x}_{k}}\right] = 4 H, \end{equation} 
and
\begin{equation}\label{eq1.8}
\left[{{x}_{k}, \left[{H, {x}_{k}}\right]}\right] = 2,
\end{equation} 
which is a version of {\it canonical commutation}.
These identities are fundamental for the subject of universal bounds 
on eigenvalues (for which see \cite{AsBe94}).   

It is argued 
in Sections~2 and 3 that from an algebraic 
point of view, the operator on a closed manifold most closely
analogous to the flat Laplacian is not the Laplace-Beltrami operator, but 
rather
\begin{equation}\label{eq1.9}
H_g := -\Delta +  g h^2 \end{equation} 
with the specific value $g = \frac14$.  

\section{Upper bounds on the fundamental gap when $M$ is a curve}

In this section, suppose that $M$ is a smooth plane or space curve
of finite length.  If $s$ denotes the 
arc length, then Schr\"odinger operators have the form
\begin{equation}\label{eq2.1}
H = -\frac{{d}^{2} }{ {ds}^{2}} + V(s).\end{equation} 
Assuming $V$ is bounded and measurable, 
$H$ is a self-adjoint operator defined 
using the Friedrichs extension from $C_c^{\infty}(M)$.  
The form domain is thereby $W_0^1(M)$, which
for a closed 
curve coincides with $W^1(M)$; otherwise
Dirichlet boundary conditions apply
at the end points.
(For standard technical facts 
about the domains of definition of 
Sturm-Liouville operators the reader is referred to \cite{ReSi80, Wei80}.)

Let $G = X_m$, the restriction to $M$ of a Euclidean coordinate 
$x_m$ in the ambient flat space, and 
calculate the first and second
commutators with $H$, recalling the classical Serret-Frenet equations:
\begin{equation}\label{eq2.2}
\left[{{H},X_{m}}\right] = -\frac{{d}^{2}X_{m} }{ {ds}^{2}} -2 \,
\frac{dX_{m} }{ ds}\frac{d }{ ds} = - \kappa {n}_{m} -2 t_{m}
\frac{d}{ ds},\end{equation} 
and
\begin{equation}\label{eq2.3}
\left[{X_{m}\left[{{H},X_{m}}\right]}\right] = 2 
{t}_{m}^{2}.\end{equation} 
There follows:

\begin{prop}\label{prop2.1} 
Let $M$ be a smooth curve in $\bbr^\nu, \nu=2$ or $3$.  
Then for $H$ as in \eqref{eq2.1} and $\varphi \in W_0^1(M),$ 
\begin{equation}\label{eq2.4}
\sum\limits_{m = 1}^{\nu } 
\left\|\left[H, X_{m}\right]\varphi \right\|^2 = 
4\int_M\left( \left|\frac{d\varphi}{ds}\right|^2
+\frac{\kappa^2 }{ 4}\left|\varphi\right|^2\right)ds.\end{equation}
\end{prop}

\begin{remn} 
\rm The assumptions on the curve can be relaxed, 
providing that the curvature and unit 
normal vector exist a.e.\ and are square-integrable.  
Likewise, potentials can be allowed 
to have certain singularities, as long as the domain of self-adjointness remains the same.
(For conditions guaranteeing this, see \cite{EdEv87, ReSi78, Wei80}.)
\end{remn}

\begin{pf} 
By closure it may be assumed that $\varphi \in C_c^{\infty}(M)$.  Apply
\eqref{eq2.2} to $\varphi$ and square the result, to obtain
$$4 \left(t_m^2\left(\frac{d\varphi }{ ds}\right)^{2}+ \frac{1}{
4}\,{\kappa }^{2}{n}_{m}^{2}\varphi^{2} + \frac{1 }{ 2}\,\kappa n_m 
t_m \varphi\,\frac{d\varphi }{ ds}\right).$$
When the sum is taken on on $m$ the last term 
drops out because {\bf t} and {\bf n} are orthogonal.  
Integration then produces Eq.~\eqref{eq2.4}.
\end{pf}

For $\varphi$ sufficiently regular, the right side of \eqref{eq2.4} equals
the one-dimensional version of
$$4 \langle{\varphi, H_{1 / 4} \varphi}\rangle.$$
Consequently, Proposition~\ref{prop2.1} implies universal 
spectral bounds for operators of the form $H_g$:

\begin{coro}\label{cor2.2} 
Let $M$ be as in {\rm Proposition~\ref{prop2.1}} and suppose that
$H$ is a Schr\"odinger Hamiltonian with $V(s)$
bounded and measurable.  Then
\begin{equation}\label{eq2.5}
\Gamma \leq 4\int_{M}{\left({{\left(\frac{d u_1  }{ ds}\right)}^{2}+ 
\frac{{\kappa }^{2} }{ 4}u_1^{2}}\right)ds}.\end{equation} 
Furthermore, if $H$ is of the form
$$H_{g} := -\frac{{d}^{2} }{ {ds}^{2}} + g {\kappa }^{2},$$
then
\begin{equation}\label{eq2.6}
\Gamma \leq \max\left({4, \frac{1 }{ g}}\right) \lambda_1. \end{equation}
 Equivalently, the universal ratio bound
$$\frac{{\lambda }_{2} }{ {\lambda }_{1}} \le 
\max\left({5,1+\frac{1 }{ g}}\right) $$
 holds.  This bound is sharp for $0 < g \leq \frac{1 }{ 4}$.
\end{coro}

\begin{remn} 
\rm In case $ g < 0,$ somewhat less simple-looking bounds on the gap in lieu 
of \eqref{eq2.6} follow easily from \eqref{eq2.5}. 
\end{remn}

\begin{pf} 
To establish \eqref{eq2.5}, make the choice 
$G = X_m$ in Lemma~\ref{lem1.1}, 
and use \eqref{eq2.3} to rewrite the left side of 
Lemma~\ref{lem1.1}.  With \eqref{eq2.4}, 
$$\Gamma \cdot 2 \sum_{m} \left\langle{{u}_{1},{t}_{m}^{2}{u}_{1}}\right\rangle
 \leq 8 \int_{M}{\left({{\left(\frac{{d{u}_{1}}  }{ ds}\right)}^{2}+ 
\frac{{\kappa }^{2} }{ 4}{{u}_{1}}^{2}}\right)ds}$$
Because $\sum_{m}{{t}_{m}^{2}} = 1$, by dividing by $2$
we arrive at  \eqref{eq2.5}.
Since 
$${\lambda }_{1} = \int_{M}\left({{\left(\frac{d{u}_{1} }{
ds}\right)}^{2}+ g{\kappa }^{2}{u}_{1}^{2} }\right)ds, $$
simple algebra yields \eqref{eq2.6}
and the equivalent ratio bound.  To see that the 
latter bound is sharp for $g \leq \frac{1 }{ 4}$, consider the case
of the circle, for which the eigenvalues are explicitly known, and:
$$\frac{{\lambda }_{2} }{ {\lambda }_{1}} =
\frac{4{\pi }^{2}\left({1+g}\right) }{
4{\pi }^{2}g} = 1+\frac{1 }{ g}.$$
\end{pf}

\section{Bounds on the fundamental gap when $M$ is a hypersurface}

This section treats the case where $M$ is a compact 
$d$-dimensional manifold smoothly immersed in $\bbr^{d+1}$.
The results are precisely as in Section~2.

\begin{thm}\label{thm3.1} 
Let $H$ be a Schr\"odinger  operator on $M$ with a bounded potential, i.e.,
\begin{equation}\label{eq3.1}
H = -\Delta + V,\end{equation} 
where $V$ is a bounded, measurable, real-valued function on $M$.  If
$M$ has a boundary, Dirichlet conditions are imposed (in the weak sense that 
$H$ is defined as the Friedrichs extension from 
$C_c^{\infty}(M)$).
Then 
\begin{align}\label{eq3.2}
\Gamma(H)  &\leq 
\frac{1 }{ d} \int_M \left({4 {\left| \nabla_{||} u_1\right|^2
+ h^2 u_1^2}}\right) dVol\nonumber\\
& = \frac{4 }{ d} \left\langle{u_1, \left({-\Delta + 
\frac{h^2 }{ 4}}\right) u_1}\right\rangle.\end{align} 
\end{thm}

\begin{remn} 
\rm The assumption of boundedness of $V$ can be considerably relaxed, as
long as the form domain of $H$ remains the Sobolev
space $W_{0}^{1}(M)$.  Conditions for this can be found in 
\cite{EdEv87, ReSi78}.
\end{remn}

\begin{pf} 
Let $\{x_m, m=0, \dots, d\}$ 
be a fixed Cartesian coordinate system for $\bbr^{d+1}$, and
denote by ${\bf e}_m$ the corresponding unit vectors.  
Ambient coordinates are introduced because they will be 
convenient for some calculations.  The 
notation $\{X_m\}$ will be used for the functions $M \to \bbr$ obtained by 
restricting $x_m$ to $M$.
Let $\nabla_{||}$ denote the tangential gradient, which can be
written in terms of the ambient gradient
by projecting away the component along ${\bf n}$:
$$\nabla_{||} = \nabla - {\bf n} \  {\bf n} \cdot  \nabla .$$
Here, {\bf n} denotes the outward normal to $M$.   
Note that the vector position on $M$, 
i.e., ${\bf X} := \sum_{m = 0}^{d} {X_m {\bf e}_m}$
is independent of the ambient coordinate system, and that
$\nabla_{||}$ is independent of the local coordinates on $M$.
Also, in an ambient 
coordinate system the Laplace-Beltrami operator $\Delta$
can be identified with $\nabla_{||}^2$.

The commutator formula \eqref{eq1.6} using $G =X_m$ states
$$ (\lambda_2 - \lambda_1) \left\langle{u_1, \left[X_m, \left[H, X_m \right] \right] 
u_1 }\right\rangle \leq 2  {\left\|\left[{H, 
{X}_{m}}\right]{u}_{1}\right\|}^{2}.$$

A simplification can be achieved after summing on $m$, and noting that 
since the potential in $H$ commutes with $X_m$, $H$ may be replaced 
in the commutators by $-\Delta$:
\begin{equation}\label{eq3.3}
 (\lambda_2 - \lambda_1)\sum\limits_{m = 0}^{d} { \left\langle{u_1, \left[X_m, 
\left[-\Delta, X_m \right] \right] u_1 }\right\rangle}  \leq 2  
\sum\limits_{m = 0}^{d} {{\left\|\left[{-\Delta, 
{X}_{m}}\right]{u}_{1}\right\|
}^{2}}. \end{equation} 

It will now be shown that the left side of 
\eqref{eq3.3} is $2 d \Gamma$.  The first step
is to realize that 
$\sum\limits_{m = 0}^{d} \left[X_m, \left[-\Delta, X_m \right] \right]$
is independent of the orientation
of the ambient Cartesian system, because 
with the orthogonality of $\{{\bf e}_m\}$ a short calculation 
shows that it is equal to the invariant expression
${\bf X} \cdot \left[{\left(-\Delta\right), {\bf X}}\right] - 
\left[{\left(-\Delta\right), {\bf X}}\right] \cdot {\bf X}$.

This allows the left side of \eqref{eq3.3} to be calculated locally in
a conveniently oriented ambient Cartesian system:  At
any point of $M$, choose the orientation
so that ${\bf n} = {\bf e}_0$, in which case
$\left[{{-\nabla }_{||}^{2}, {X}_{m}}\right] = -2 \frac{\partial }{ {\partial 
x_m}}+\left({\rm scalar \ function}\right)$, $m=1, \dots, d$, while 
$\left[{{-\nabla }_{||}^{2}, {X}_{0}}\right]$ is a scalar function.
Consequently, 
$$\sum\limits_{m = 0}^{d} {\left[X_m, 
\left[{-\Delta,X_{m}}\right]\right]} = 
\sum\limits_{m = 0}^{d} {\left[X_m, 
\left[{{-\nabla }_{||}^{2},{X}_{m}}\right]\right]} = 2d.$$
By integrating, the left side of \eqref{eq3.3} is therefore
$2 d \left(\lambda_2 - \lambda_1\right) = 2 d \Gamma$.

To facilitate the calculation of the right side of
\eqref{eq3.3}, define an operator 
${\bf P}$ mapping the domain of $H$ to $\bbr^{d+1} \otimes L^2(M)$ via 
$${\bf P} f:=  
- \frac{1 }{ 2} \sum\limits_{m = 0}^{d} {\bf e}_m \left[-\Delta, 
X_m\right]  f =  \frac{1 }{ 2}  \left[\Delta, {\bf X}\right]  f.$$
Observe now that  due to the orthonormality of 
$\{ {\bf e}_m \}$, the right side of \eqref{eq3.3} is 
$$4 \left\|{\bf P} u_1\right\|_{\bbr^{d+1} 
\otimes L^2(M)}^{2}.$$
The Laplace-Beltrami operator can be conveniently expressed 
in local Riemann normal coordinates on $M$, in 
which it is arranged that at a chosen point $p \in M$, 
$\Delta = \sum_{j = 1}^{d} \frac{\partial^2 }{\partial s_j},$
where $\{s_j\}$ measure arc length along the lines of curvature passing through 
$p$.  Thus, at $p$,
\begin{equation}\label{eq3.4}
\left[-\Delta, X_m\right] = \sum\limits_{j = 1}^{d}{\left(-2  
\frac{\partial X_m }{ \partial s_j}
\frac{\partial }{ \partial s_j} - 
\frac{\partial^2 X_m }{ \partial s_j^2}\right)},\end{equation} 
so
\begin{align}\label{eq3.5}
{\bf P} u_1&=   \sum\limits_{j = 1}^{d}{\left(
\frac{\partial {\bf X} }{ \partial 
s_j}\frac{\partial u_1 }{ \partial s_j} + \frac{1 }{ 2}
\frac{\partial^2 {\bf X} }{ \partial 
s_j^2}u_1\right)}\nonumber\\
&= \sum\limits_{j = 1}^{d}{\left({{\bf t}_j}
\frac{\partial u_1 }{ \partial s_j} 
\pm \frac{1 }{ 2} \kappa_j {\bf n} u_1\right).}\end{align} 

The second line of \eqref{eq3.5} uses the
Serret-Frenet equations, according to which ${\bf t}_j$ is a unit tangent 
vector parallel to the 
line of curvature parametrized by ${s_j}$, and
$\kappa_j$ is the associated principal curvature.  
In other words, $\{{\bf t}_j\}$ are the normalized eigenvectors of the 
shape operator.  The varying sign enters
because here ${\bf n}$ is 
by convention outward, whereas the normal defined by the Serret-Frenet 
equations may be outward or inward.   
Where the latter is not defined, its coefficient is 0,
and it drops out.

Since $\{ {\bf n}, {\bf t}_j\}$ is an orthonormal system, at the point $p$,
$$\left\|{\bf P} u_1\right\|_{\bbr^{d+1}}^{2}\left(p\right) = 
 \sum\limits_{j = 1}^{d}\frac{\partial u_1 }{ \partial s_j^2}
+ \frac{1 }{ 4}  \left(\sum\limits_{j = 1}^{d} {\kappa_j} u_1\right)^2.$$
This formula is equal to the coordinate-independent expression
\begin{equation}\label{eq3.6}
\left| \nabla_{||} u_1\right|^2+ \frac{h^2 }{ 4} u_1^2,\end{equation}
which can now be integrated over $M$.  
When \eqref{eq3.6} is used in the right side of \eqref{eq3.3}, 
which is then divided by $2 d$, the result is
$$(\lambda_2 - \lambda_1) \leq \frac{4 }{ d} 
\int_M \left({{\left| \nabla_{||} u_1\right|^2
+ \frac{h^2 }{ 4} u_1^2}}\right) dVol.$$
This is equivalent to \eqref{eq3.2}.
\end{pf}

 \smallskip
An immediate consequence of Theorem~\ref{thm3.1} by partial integration is:

\begin{cor}\label{cor3.2} 
Let $H$ be as in \eqref{eq3.1} and define 
$\delta := \sup_M \left({\frac{h^2 }{ 4} - V}\right)$.
Then
$$\Gamma(H) \leq \frac{4 }{ d} \left({\lambda_1 + \delta}\right).$$
\end{cor}

A further corollary is an
isoperimetric spectral theorem for operators of the form
$H_g$ from \eqref{eq1.9}:

\begin{cor}\label{cor3.3} 
Let $H_g$ be defined on $M$, a $d$-dimensional manifold smoothly 
immersed in $\bbr^{d+1}$. 
Then for $0 < g \leq \frac{1 }{ 4}$, the eigenvalues satisfy
\begin{equation}\label{eq3.7}
{\lambda_2 - \lambda_1} \leq \frac{4 \sigma \lambda_1 }{ d}\,,\end{equation} 
with
\begin{equation}\label{eq3.8}
\sigma = \max\left({1, \frac{1 }{ {4 g}}}\right).\end{equation}
For the ratio, these bounds read
$$\frac{\lambda_2 }{ \lambda_1} 
\leq \frac{{1 + gd} }{ gd}\, , \quad\mbox{for }  0 < g \leq \frac{1 }{ 4} $$
and
$$\frac{\lambda_2 }{ \lambda_1} \leq \frac{{4 + d} }{ d}\,,  
\quad\mbox{for }  g \geq \frac{1 }{ 4}. $$
The bounds \eqref{eq3.7} are optimal for $g \leq \frac{1 }{ 4}.$
\end{cor}

\begin{pf} 
If $H = H_g$, then
$$\int_M{\left({\left| \nabla_{||} u_1\right|^2
+ {g h^2} u_1^2}\right) dVol} = \left(u_1, H_g u_1\right) = \lambda_1,$$
a multiple of which majorizes the right side of \eqref{eq3.2}:  When
$g \leq \frac{ 1 }{ 4}$, this is immediate.
Otherwise first multiply 
$\left| \nabla_{||} u_1\right|^2$ by $\frac{1 }{ 4g} > 1$.  The result is
\eqref{eq3.7}.

Optimality is established for $g \leq \frac{1 }{ 4}$ 
by the example of the sphere $S^d,$ for which, 
after normalizing  the scale so that $S^d$ embeds in 
$\bbr^{d+1}$ as the unit sphere,
$$\lambda_1 = g d^2,\quad  \lambda_2 = g d^2 + d$$
\cite{Mul66}.  Thus 
$$d = \lambda_2 - \lambda_1 \leq  \left(\frac{gd^2}{
gd}\right) = d.$$
\end{pf} 

\section{Sum rules for eigenvalues}

In \cite{HaSt97} universal estimates for Laplacians and Schr\"odinger 
operators were obtained by deriving {\it sum rules} with commutator relations,
and appropriately modified versions of these apply to operators of the type 
\eqref{eq3.1}.
This section follows that article closely, but with the purpose of
displaying the effect of mean curvature.  
The results will subsume those of Section~3 and give information
about the spectrum above the first two eigenvalues.

For $H$ as in \eqref{eq3.1} the calculations leading to 
\eqref{eq3.5} and \eqref{eq3.6}
established the operator identities	 
\begin{equation}\label{eq4.1}
{\bf P} = \sum\limits_{j = 1}^{d}{\left({{\bf t}_j}
\frac{\partial  }{ \partial s_j} 
\pm \frac{1 }{ 2} \kappa_j {\bf n} \right)}\end{equation} 
and for a dense set of functions $\varphi$, 
\begin{equation}\label{eq4.2}
\|{\bf P} \varphi\|^2 = 
\left\langle{\varphi, H_{1 / 4} \varphi}\right\rangle.\end{equation} 
Thus ${\bf P}$ plays the r\^ole of a
momentum operator, with which
there is a version of canonical commutation (cf.~\eqref{eq1.8})
as follows.  Defining a variant commutator bracket for
operators ${L^2(M)} \to \bbr^{d+1} \otimes L^2(M)$ by 
$\left[{A; B}\right] := A \cdot B - B \cdot A,$
a calculation shows that
 $\left[{{\bf P}; X_k {\bf e}_k}\right] = \sum_{j=1}^d{{\bf t}_j \cdot
\frac{{\partial  X_k {\bf e}_k} }{ {\partial s_j}}} = {\bf 1}$
(identity operator), and by averaging on $k$,
\begin{equation}\label{eq4.3}
 {\bf 1} = \frac{1 }{ d}\,\left[{{\bf P}; {\bf X}}\right].\end{equation}
This is a coordinate-independent formula.

A sum rule for $H$, analogous to that of \cite{HaSt97} reads as follows:

\begin{prop}\label{prop4.1} 
Let $H$ be as in \eqref{eq3.1}, with eigenvalues $\{{\lambda_k}\}$
and normalized eigenfunctions $\{ {u_k} \}$.  Then
\begin{equation}\label{eq4.4}
1 = \frac{4 }{ d}
\sum_{\stackrel{\scriptstyle k}{\lambda _{k}\ne \lambda_j}}
\frac{\left|\left\langle u_k,{\bf P}u_j\right\rangle\right|^{2} }
{\lambda_k-\lambda_j}.\end{equation} 
Furthermore, if $f$ is any function summable on the spectrum
$\sigma(H),$ then
\begin{equation}\label{eq4.5}
\sum\limits_{j}^{\infty } f\left({{\lambda }_{j}}\right) =
-\frac{2 }{ d}\, \sum\limits_{\stackrel{\scriptstyle j,k}
{\lambda_k\ne \lambda_j}} 
{{\left|{\left\langle{{u}_{k},{\bf P}{u}_{j}}\right\rangle}\right| }^{2}\,
\frac{{ f\left({\lambda_j}\right) - f\left({\lambda_k}\right)}
 }{ {{\lambda }_{j}-{\lambda }_{k}}}}.\end{equation} 
\end{prop}

\begin{pf} 
   From the commutator gap formula \eqref{eq1.2}, with $G = {\bf X}$, 
$$- 2 \left\langle{u_j, {\bf P} u_k}\right\rangle = (\lambda_j -\lambda_k) \left\langle{u_j, {\bf X} u_k}\right\rangle,$$
so for $\lambda_j \ne \lambda_k,$
\begin{equation}\label{eq4.6}
\left\langle{{u}_j , {\bf X}{u}_{k}}\right\rangle = 
-2\frac{\left\langle{{u}_{j},{\bf P} u_{k}}\right\rangle }{
{\lambda }_{j}-{\lambda }_{k}}.\end{equation} 
Meanwhile, from \eqref{eq4.3} and the completeness of the
orthonormal set of real-valued eigenfunctions
$\left\{{u_k}\right\}$, it follows that
\begin{align*}
1 &= \left\langle{{u}_{j},{u}_{j}}\right\rangle = \frac{1 }{ d} 
\left\langle{{u}_{j},\left[{{\bf P}; 
{\bf X}}\right]{u}_{j}}\right\rangle
= \frac{2 }{ d} \left\langle{{u}_{j},{\bf P} 
\cdot {\bf X} {u}_{j}}\right\rangle\\
&= \frac{2 }{ d} \sum_{k}{  \left\langle{u_j, {\bf P} u_k}\right\rangle
\cdot \left\langle{u_k, {\bf X} u_j}\right\rangle}.
\end{align*}
Substitution from \eqref{eq4.6} yields \eqref{eq4.4}.  
Observe that there is no sum on $j$ in \eqref{eq4.4}; 
it is true for all $j.$
To derive the identity \eqref{eq4.5} from \eqref{eq4.4}, sum on $j$ to obtain
$$\sum\limits_{j}^{\infty } f\left({{\lambda }_{j}}\right) =
-\frac{4 }{ d}\sum\limits_{\stackrel{\scriptstyle j,k}
{\lambda_k\ne \lambda_j}} 
{{\left|{\left\langle{{u}_{k},{\bf P}{u}_{j}}\right\rangle}\right| }^{2}\,
\frac{{ f\left({\lambda_j}\right)}
 }{
{{\lambda }_{j}-{\lambda }_{k}}}}, $$
and then symmetrize in the indices $j$ and $k$.
\end{pf}

\smallskip
Because \eqref{eq4.4} and \eqref{eq4.5} are identical in form  
to the sum rules derived in \cite{HaSt97}, 
the corollaries in that
article can be carried over directly.  
Accordingly, the proofs of their analogues
will be given in outline only.
For operators 
$H_g$, the
constant $\sigma$ in Eq.~(10) of \cite{HaSt97} is
$\sigma = \max\left({1, \frac{1 }{ {4 g}}}\right)$
as in \eqref{eq3.8}.  
Corollary~3 of \cite{HaSt97} is a kind of
Hile-Protter inequality \cite{HiPr80}.  Its analogue 
is as follows:

\begin{coro}\label{cor4.2} 
Let $H$ be as in \eqref{eq3.1} and 
$\delta := \sup_M \left({\frac{h^2 }{ 4} - V}\right)$ 
as in {\rm Corollary~\ref{cor3.2}}.  Then for each $n = 1,2, \dots$,
such that $\lambda_{n+1} \ne \lambda_n$,
\begin{equation}\label{eq4.7}
1 \leq \frac{4 }{ {d n}} \sum_{j=1}^{n}{}
\frac{{\lambda_j + \delta} }{ {\lambda_{n+1} - \lambda_j}}.\end{equation}
Suppose that $H_g$ is as in \eqref{eq1.9} with $M$ as in 
{\rm Section~3} and $g > 0$.  If $\lambda_{n+1} \ne \lambda_{n}$, then
\begin{equation}\label{eq4.8}
1 \leq \frac{4 \sigma  }{ {d n}} 
\sum_{j=1}^{n}{}\frac{{\lambda_j} }{ {\lambda_{n+1} - \lambda_j}}.
\end{equation} 
\end{coro}

\begin{pf} 
Sum \eqref{eq4.4} on $j$ from $1$ to $n$.  Since by symmetry, 
$$\sum\limits_{\stackrel{\scriptstyle j,k \leq n}
{\lambda_k\ne \lambda_j}} \frac{ 
{\left|{\left\langle{{u}_{k},{\bf P}{u}_{j}}\right\rangle}\right| }^{2} }{
{\lambda }_{k}-{\lambda }_{j}} = 0,$$
this produces
$$
n = 
\frac{4 }{ d} \sum\limits_{\stackrel{\scriptstyle j \leq n}
{k \geq n+1}} \frac{ 
{\left|{\left\langle{{u}_{k},{\bf P}{u}_{j}}\right\rangle}\right| }^{2} }{
{\lambda}_{k}-{\lambda}_{j}},$$
so
$$1 \leq 
\sum\limits_{j = 1}^{n} \frac{ 
1 }{
{{\lambda}_{n+1}-{\lambda}_{j}}}  \sum\limits_{k=0}^{\infty}
{\left|{\left\langle{{u}_{k},{\bf P}{u}_{j}}\right\rangle}\right|^2} = 
\sum\limits_{j = 1}^{n} \frac{{\|{\bf P}{u}_{j}\|^2}
}{
{\lambda_{n+1}-{\lambda}_{j}}}.
$$
Inequality \eqref{eq4.7} follows by substituting from \eqref{eq4.2}
 and using the elementary fact that 
\begin{equation}\label{eq4.9}
H_{1 / 4} = H  + \frac{h^2 }{ 4} - V \leq H + \delta
\end{equation}
in the sense of quadratic forms.

For the universal inequality \eqref{eq4.8}, we instead use
the fact that $g' \geq g$ implies $H_{g'} \geq H_g$
and, for $g < \frac{1 }{ 4}$, that
$H_{1 / 4} \leq \frac{1 }{ {4 g}} H_g$, again in the sense
of quadratic forms.
\end{pf}

\smallskip
Theorem~5 of \cite{HaSt97} is a ``Yang-type"
inequality \cite{Yan91}, which in the 
context of this article reads:

 \begin{coro}
 \label{cor4.4} 
 Let $H$ be as in \eqref{eq3.1}, with $M$ a compact, smooth submanifold.
 For $z$ such that $\lambda_n < z \leq \lambda_n+1$,
 \begin{align*}
 \sum\limits_{j =1}^{n}{(z -\lambda_j)^{2}}&\le
 \frac{4 }{ {n d}}\sum\limits_{j= 1}^{n}
 {(z -\lambda_j){\|{\bf P} u_j\|}^2}\\
 &\leq  \frac{4 }{ {n d}}\sum\limits_{j= 1}^{n}
 (z -\lambda_j)(\lambda_j+\delta).\end{align*}
 \end{coro}

\begin{pf} 
This results from multiplying \eqref{eq4.4} by
$(z -\lambda_j)^{2}$ and summing on $j$ from 1 to $n$.   The sum on
$k$ from 1 to $n$ can be dropped because the summand is antisymmetric 
under exchange of $j$ and $k$.   For $\lambda_k \geq z$, 
$\frac{\left({z - \lambda_j}\right)^{2}} { \lambda_k - \lambda_j} 
\leq \left({z - \lambda_j}\right),$
whence
\begin{align*}
\sum\limits_{j =1}^{n} (z -\lambda_j)^{2}&\leq 
\frac{4 }{ {n d}}\sum\limits_{j= 1}^{n} 
{(z -\lambda_j)\sum\limits_{k=n+1}^{\infty} 
{\left|\left\langle{u_k, {\bf P} u_j}\right\rangle\right|}^2}\\
&\leq \frac{4 }{ {n d}}\sum\limits_{j= 1}^{n}
{(z -\lambda_j){\|{\bf P} u_j\|}^2}\end{align*}
by Bessel's inequality.  The final inequality
follows with \eqref{eq4.9}
\end{pf} 

\smallskip
Since Corollary~\ref{cor4.4}
 states that a certain quadratic function of $z$ is negative when 
$\lambda_n < z \leq \lambda_n+1$, 
it implies that each spectral gap can be bounded by
an explicit expression in terms of the distribution of lower eigenvalues.  
Define
$$\overline{\lambda_n} := \frac{1 }{ n} \sum\limits_{k=1}^{n} {\lambda_k}$$
and
$$\overline{\lambda_{n}^{2}} := 
\frac{1 }{ n} \sum\limits_{k=1}^{n} {\lambda_{k}^2}.$$
Then each eigenvalue gap of a generic $H$ is bounded by a complicated, explicit 
expression where the mean curvature enters through $\delta$, whereas the 
eigenvalue gaps of operators of the form $H_g$ are bounded by a 
simpler universal expression:

\begin{coro}\label{cor4.5} 
\smallskip
\begin{itemize}
\item[\rm a)]  Let $H$ be as 
\eqref{eq3.1}, with $M$ a compact, smooth submanifold.  Then
$$\left[{\lambda_n, \lambda_{n+1}}\right] 
\subseteq \left[{\left({1+\frac{2 }{ d}}\right)\overline{{\lambda }_{n}}-\sqrt{D_n^{\delta}}, \left({1+\frac{2 }{ d}}\right)\overline{{\lambda }_{n}}+ \sqrt{D_n^{\delta}}}\right],$$
with
$$D_n^{\delta} := 
 \frac{4 }{ {d^2}}\left({{\left({\frac{dn+2 }{ 2}}\right)}^{2}
{\overline{{\lambda }_{n}}}^{2} + \left({dn-d+2}\right)
\delta \overline{{\lambda }_{n}}
-d\left({\frac{dn+4 }{ 4}}\right)\overline{{\lambda }_{n}^{2}}
 + {\delta }^{2}}\right).$$
\item[\rm b)]  For $H_g$, $0 < g$, of the form \eqref{eq1.9} on a 
smooth, compact submanifold $M$, 
$$\left[{\lambda_n, \lambda_{n+1}}\right]
\subseteq \left[{\left({1+\frac{2 \sigma}{ d}}\right)\overline{{\lambda }_{n}}-\sqrt{D_n}, 
\left({1+\frac{2 \sigma}{ d}}\right)\overline{{\lambda }_{n}}+ \sqrt{D_n}}\right],$$
with
$$D_n:= \left(\left(1+\frac{2\sigma  }{ d}\right)\overline{\lambda_n}\right)^{2} - 
\left(1+\frac{4\sigma  }{ d}\right)\overline{\lambda_{n}^{2}} > 0. $$
\noindent
This bound is sharp for every non--zero eigenvalue gap of
$H_{\frac{1}{4}}$ on the sphere.
\end{itemize}
\end{coro}

\begin{pf}  
It is only necessary to calculate the roots of the quadratic 
expression from the preceding Corollary,
$$z \rightarrow \sum\limits_{j =1}^{n} (z -\lambda_j)^{2} - \frac{4 }{ {n d}}\sum\limits_{j= 1}^{n}
{(z -\lambda_j){\|{\bf P} u_j\|}^2}$$
and to substitute from \eqref{eq4.2}.
(Cf.~\cite{HaSt97}, Proposition~6;  the bound b) 
is in fact identical in form to the one in that article.)

An explicit calculation shows
that the bound is sharp for the non-zero eigenvalue gaps of the sphere,
for which all the eigenvalues are known and elementary \cite{Mul66}:
For simplicity, assume that $d = 2, g = \frac{1}{4}$, and that $M$ is
the sphere of radius 1 
embedded in $\bbr^{3}$.  Then $h = 2, \sigma = 1$, and:

$$\lambda_1 = 1; \lambda_2 = \lambda_3 = \lambda_4 = 3; \dots; \lambda_{\left({m-1}\right)^2+1} = 
\dots = \lambda_{m^2} = m^2 - m + 1.$$   
\noindent
For $n = m^2$, the calculation shows that 
$\overline{\lambda_n} = \frac{n+1}{2}$, and $\overline{{\lambda }_{n}^{2}} = \frac{n^2+n+1}{3}.$
\noindent
Hence $D_n = n$, and b) informs us that 

$$2 \overline{\lambda_{m^2}} - m =  m^2 - m + 1 \leq  \lambda_{m^2} =  m^2 - m + 1$$
$$\ \ \ \ \ \ \ \   \leq  \lambda_{m^2+1} = m^2 + m + 1 \leq  2 \overline{\lambda_{m^2}} + m = m^2 + m + 1,$$

\noindent
and thus $\lambda_{m^2}$ equals the lower bound $2 \overline{\lambda_{m^2}} - m$ and $\lambda_{m^2+1}$
equals the upper bound $2 \overline{\lambda_{m^2}} + m$.

\end{pf}

\smallskip
Finally, consider the partition function for $H$,
$$        Z(t) := \text{tr}(\exp(-tH)),                                                 
$$
If the function $f$ of Proposition~{4.1} is chosen as $f(x) := \exp(-tx)$, 
then (after a short calculation exactly as for Eq. (15) of 
\cite{HaSt97}):

\begin{equation}\label{eq4.10}
{Z\left(t\right) \le  \left({\frac{2 t}{d}}\right) \sum\limits_{
j}  \left(\exp \left(-t \lambda_j\right)\right) 
){\|{\bf P} u_j\|}^2},
\end{equation}

\noindent
which implies the following bounds:

\begin{coro}\label{cor4.6} 
\begin{itemize}
\item[\rm a)]  Let $H$ be as in
\eqref{eq3.1}, with $M$ a compact, smooth submanifold.  Then
$t^{\frac{d}{2}} \exp \left({- \delta t}\right) Z(t)$ is a nonincreasing function;
\item[\rm b)]  For $H_g$ be of the form \eqref{eq1.9} on a 
smooth, compact submanifold $M$, 
$t^{\frac{d}{2 \sigma}} Z(t)$ is a nonincreasing function.
\end{itemize}
\end{coro}

\begin{pf}  
For general $H$, the right side of \eqref{eq4.10} is majorized by 
$$
\left({\frac{2 t}{d}}\right) \sum\limits_{
j}{\lambda_j\exp \left(-t \lambda_j\right)} + \left({\frac{2 t \delta}{d}}\right) Z(t).$$
Hence
$$
 \left({1 - {\frac{2 t \delta}{d}}}\right) Z(t)  
+\left({\frac{2 t}{d}}\right) Z^{\prime}(t) \leq 0,$$
and when multiplied by the integrating factor 
$t^{\frac{d}{2}} \exp \left({- \delta t}\right)$
the left side becomes a positive multiple of the derivative of
$t^{\frac{d}{2}} \exp \left({- \delta t}\right) Z(t)$.

Statement {\rm b)} follows similarly, after majorizing the 
right side of  \eqref{eq4.10} with 
$$
\left({\frac{2 \sigma t}{d}}\right) \sum\limits_{
j}  \left(\exp \left(-t \lambda_j\right)\right) 
)\lambda_j.$$
\end{pf}


\paragraph*{Acknowledgments}
{\rm Part of this work was done while the author visited CEREMATH, at the
Universit\'e de Toulouse 1.  NSF support through grant DMS-0204059 is 
also gratefully acknowledged.}

\end{document}